\documentclass[fleqn]{mat01}
\usepackage{times,mathtimy,amssymb}
\begin{document}

\font\xx=msam5 at 9pt
\def\ab{\mbox{\xx{\char'03}}}
\def\ba{\mbox{\xx{\char'04}}}

\font\bi=tibi at 10.4pt

\setcounter{page}{71}
\firstpage{71}


\def\ques{\trivlist\item[\hskip\labelsep{\it Question.}]}
\def\defi{\trivlist\item[\hskip\labelsep{\bf DEFINITION.}]}
\def\remark{\trivlist\item[\hskip\labelsep{\it Remark.}]}
\def\remarks{\trivlist\item[\hskip\labelsep{\it Remarks}]}
\def\noot{\trivlist\item[\hskip\labelsep{{\it Note.}}]}
\def\thoe{\trivlist\item[\hskip\labelsep{{\bf Theorem}}]}
\newtheorem{theo}{Theorem}
\renewcommand\thetheo{\arabic{section}.\arabic{theo}}
\newtheorem{theor}{\bf Theorem}
\newtheorem{lem}{Lemma}
\newtheorem{fact}{Fact}
\newtheorem{definit}{\rm DEFINITION}
\newtheorem{propo}{\rm PROPOSITION}
\newtheorem{rema}{Remark}
\newtheorem{exam}{Example}
\newtheorem{coro}{\rm COROLLARY}

\title{Questions concerning matrix algebras and
invariance\\ of~spectrum}

\markboth{Bruce A Barnes}{Questions concerning matrix algebras and
invariance of spectrum}

\author{BRUCE A BARNES}

\address{Department of Mathematics, University of Oregon, Eugene, Oregon 97403, USA\\
\noindent E-mail: barnes@darkwing.uoregon.edu}

\volume{113}

\mon{February}

\parts{1}

\Date{MS received 2 January 2002; revised 2 April 2002}

\begin{abstract}
Let $A$ and $B$ be unital Banach algebras with $A$ a subalgebra of $B$.
Denote the algebra of all $n\times n$ matrices with entries from $A$ by
$M_{n}(A)$. In this paper we prove some results concerning the open
question: If $A$ is inverse closed in $B$, then is $M_{n}(A)$ inverse
closed in $M_{n}(B)$? We also study related questions in the setting
where $A$ is a symmetric Banach *-algebra.
\end{abstract}

\keyword{Banach algebra; inverse closed; symmetric *-algebra; matrix
algebra.}

\maketitle

\section{Introduction}

Throughout, $A$ and $B$ are unital Banach algebras, $A$ is a subalgebra
of $B$,   and the unit of $B$ is in $A.$ In this paper, we continue
studying the relationships among the three properties:
\renewcommand\thefootnote{}\footnote{Dedicated to Professor Ashoke K.
Roy on his retirement.} 

{\it $A$ is inverse closed in $B$} if whenever $a\in A $ and $a^{-1}\in
B$, then $a^{-1}\in A$.

{\it $A$ is SRP in $B$} if $r_{A}(a)=r_{B}(a)$ for all $a\in A$ \ (here
$r_{A}(a)$ is the spectral radius of $a\in A$ in $A;$ SRP stands for
`spectral radius preserving').

When $A$ has an involution $\ast $, {\it $A$ is *-inverse closed in $B$}
if whenever $a=a^{\ast }\in A$ and $a^{-1}\in B$,   then $a^{-1}\in A$.

Note that $A$ is inverse closed in $B$ is equivalent to: \ for all $a\in
A$,   $\sigma (a;A)=\sigma (a;B)$ (here $\sigma (a;A)$ denotes the
spectrum of $a$ relative to the algebra $A)$.

A Banach *-algebra $A$ with unit $1$ is {\it symmetric} if $
1+a^{\ast }a$ is invertible \ for all $a\in A$. This is equivalent to
the property that $\sigma (a^{\ast }a;A)\subseteq \lbrack 0,\infty )$
for all $ a\in A$. The relationships among the three properties above
when $A$ is a symmetric Banach *-algebra is the subject of the author's
paper \cite{[B2]}. In this paper we study these properties in algebras of
$n\times n$ matrices over a Banach algebra; notation: $M_{n}(A)$ denotes
the algebra of all $n\times n$ matrices with entries from $A$. Much of
our work in this paper centers on the open question:

\begin{ques}
If $A$ is inverse closed in $B$, then is
$M_{n}(A)$ inverse closed in $M_{n}(B)$?
\end{ques}\vspace{.3pc}

It is known that when $A$ is commutative, then this question has an
affirmative answer; this follows from [6, Theorem~1.1]. Many other results
related to this question are known. We list two useful results; Fact 1
is proved in [9, Theorem~2.2.14, p. 219], and Fact 2 is proved in [3,
Proposition~2]. 

\begin{fact}
{\rm If $A$ is SRP in $B$, then $M_{n}(A)$ is SRP in $M_{n}(B)$.}
\end{fact}

\begin{fact}
{\rm Assume that $A$ is SRP in $B$. Let $\overline{A}$ be the closure
of} $A$ i{\rm n} $B$. {\rm If} $a\in A$ {\rm and} $a^{-1}\in \overline{A}$, {\rm then} $a^{-1}\in
A$.
\end{fact}

The matrix algebras $M_{n}(A)$ and $M_{n}(B)$ are Banach algebras with
respect to natural norms: For $(D,\Vert $ $\Vert _{D})$ a Banach
algebra, define for $T=(t_{jk})\in M_{n}(D)$,   
\begin{equation*}
\Vert T\Vert =\sum_{j=1}^{n}\sum_{k=1}^{n}\Vert t_{jk}\Vert _{D}.
\end{equation*}
Then $(M_{n}(D),\Vert $ $\Vert )$ is a Banach algebra. We use the norm
defined above, or any equivalent norm, as the standard norm on
$M_{n}(D)$.

Combining Facts 1 and 2, we have an affirmative answer to the Question
when $A$ is dense in $B$. (The proof is easy: Assume that $A$ is
inverse closed in $B$,   and $A$ is dense in $B$. By Fact 1, $M_{n}(A)$ is
SRP in $M_{n}(B)$. Since $A$ is dense in $B$, it is easy to see that
$M_{n}(A)$ is dense in $M_{n}(B)$. It follows from Fact 2 that
$M_{n}(A)$ is inverse closed in $M_{n}(B))$. We state this result as a
proposition.

\setcounter{propo}{2}
\begin{propo}$\left.\right.$\vspace{.3pc}

\noindent If $A$ is inverse closed in $B$,   and $A$ is dense in $B$, then
$M_{n}(A)$ is  inverse closed in $M_{n}(B)$.
\end{propo}

\section{Results for a general Banach algebra {\bi A}}

For a Banach algebra $D$, we let $\hbox{Inv}_{\rm l}(D), \hbox{Inv}_{\rm r}(D)$, and
$\hbox{Inv}(D)$ denote the set of left invertible elements of $D$, the set of
right invertible elements of $D$, and the set of invertible elements of
$D$,   respectively. Also, we write $GL_{n}(A)=\hbox{Inv}\ (M_{n}(A))$. 

We introduce an equivalence relation on $M_{n}(B)$: For $T,S\in
M_{n}(B), $ we write $T\thickapprox S(GL_{n}(A))$ if there exist
$V,W\in GL_{n}(A)$ such that $VTW=S$. It is clear that this is an
equivalence relation on $M_{n}(B)$. For convenience we usually write
$T\thickapprox S$ with the expression $GL_{n}(A)$ omitted. Since
$GL_{n}(A)$ is a group, if $T\in M_{n}(A),\;S\in GL_{n}(A)$,   and
$T\thickapprox S$, then $T\in GL_{n}(A). $ In particular, it is easy to
check that if $T\in M_{n}(A)$ and $S$ is obtained from $T$ by a
finite sequence of interchanges of two rows or two columns, then
$T\thickapprox S$. For example,
\begin{equation*}\left( 
\begin{array}{cc}
a & b \\ 
c & d
\end{array}
\right) \thickapprox \left( 
\begin{array}{cc}
c & d \\ 
a & b
\end{array}
\right),
\end{equation*}
since
\begin{equation*}
\left( 
\begin{array}{cc}
a & b \\ 
c & d
\end{array}
\right) =\left( 
\begin{array}{cc}
0 & 1 \\ 
1 & 0
\end{array}
\right) \left( 
\begin{array}{cc}
c & d \\ 
a & b
\end{array}
\right).
\end{equation*}

Now we prove a preliminary result which we believe is known (we have
been told that this proposition follows from results in \cite{[KM]}).

\begin{propo}$\left.\right.$\vspace{.3pc}

\noindent Assume that $A$ is inverse closed in $B${\rm ,} and $A$ is continuously
embedded in$\;B$. {\rm (}`$A$ is continuously embedded in $B$' means that \
there exists $J>0$ \ such that $J$ $\Vert a\Vert_{A}$ $\ \geq \Vert
a\Vert _{B}$ for all $a\in A$.{\rm )} Let $T=(t_{jk})\in M_{2}(A)$ with 
$T^{-1}\in $ $M_{2}(B)$. Also assume that some entry in $T$ is contained
in $\overline{{\hbox{Inv}_{\rm l}(A)}}\cup \overline{{\hbox{Inv}_{\rm r}(A)}}$. Then $T^{-1}\in
M_{2}(A)$. 
\end{propo}

\begin{proof}
By interchanging rows and columns of $T$ if necessary, we may assume
that $t_{11}\in {\hbox{Inv}_{\rm l}(A)}$. First assume that $t_{11}\in
\hbox{Inv}_{\rm l}(A)$. Choose $a\in A$ with $at_{11}=1$, and choose $\lambda \in
\mathbf{C}$ such that $\lambda +t_{21}\in \hbox{Inv}(A)$. Let 
\pagebreak

$\left.\right.$\vspace{-.8cm}

\begin{equation*}
R=\left( 
\begin{array}{cc}
\lambda a & 1 \\ 
1 & 0
\end{array}
\right) T=\left( 
\begin{array}{cc}
\lambda +t_{21} & r_{12} \\ 
t_{11} & t_{12}
\end{array}
\right).
\end{equation*}
Note that by construction, $r_{11}=\lambda +t_{21}\in \hbox{Inv}(A)$. Let
\begin{equation*}
S=\left( 
\begin{array}{cc}
1 & 0 \\ 
-t_{11}r_{11}^{-1} & 1
\end{array}
\right) R=\left( 
\begin{array}{cc}
\lambda +t_{21} & r_{12} \\ 
0 & s_{22}
\end{array}
\right).
\end{equation*}
Note that $V=\left( 
\begin{array}{cc}
\lambda a & 1 \\ 
1 & 0
\end{array}
\right) $ and $W=\left( 
\begin{array}{cc}
1 & 0 \\ 
-t_{11}r_{11}^{-1} & 1
\end{array}
\right) $ are in $GL_{n}(A)$. Also, $S=WVT$.

Thus, $S^{-1}\in M_{2}(B)$, and $\lambda +t_{21}\in \hbox{Inv}(B)$. It
follows from [9, the criterion on p.~78] that $s_{22}\in
\hbox{Inv}(B)$. Since $A$ is inverse closed in $B$,   $\lambda
+t_{21},\;s_{22}\in \hbox{Inv}(A)$. It follows from a straightforward
computation that $S\in GL_{n}(A)$. Also, $T\thickapprox S$,   and
therefore, $T^{-1}\in M_{2}(A)$. 

Now assume that $t_{11}\in \overline{{\hbox{Inv}_{\rm l}(A)}}.\;$Choose $
\{a_{m}\}\subseteq \hbox{Inv}_{\rm l}(A)$ with $\Vert t_{11}-a_{m}\Vert
_{A}\rightarrow 0$, and so $\Vert t_{11}-a_{m}\Vert _{B}\rightarrow 0$.
Since $GL_{n}(B)$ is open, we may assume for all $m$,
\begin{equation*}
T_{m}=\left( \begin{array}{cc} a_{m} & t_{12} \\
 t_{21} & t_{22}
\end{array}\right) \in GL_{n}(B).
\end{equation*}
By the previous argument, $T_{m}$
$\in $ $GL_{n}(A)$ for all $m$. Also, $T_{m}^{-1}\rightarrow T^{-1}$ in
$M_{2}(B)$. It follows from Facts 1 and 2 that $T^{-1}\in M_{2}(A)$.
\hfill \ba
\end{proof}

Let $U_{n}(A)$ be the algebra of upper triangular matrices in 
$M_{n}(A)$,   that is, $T=(t_{jk})\in U_{n}(A)$ if $t_{jk}=0$ whenever
$j>k$. When $A$ is inverse closed in $B$, it is easy to check that
$U_{n}(A)$ is inverse closed in $U_{n}(B)$. For example, suppose in the
case $n=2$ \ that 
$T=\left( \begin{array}{cc} 
t_{11} & t_{12} \\
 0 &t_{22}
\end{array}\right)
 \in U_{2}(A)$ has inverse 
$S=\left(
\begin{array}{cc}
s_{11} & s_{12} \\
 0 & s_{22}
 \end{array} \right)
\in U_{2}(B)$. 
Then $ST=\left( \begin{array}{cc} 
1 & 0 \\
 0 & 1
\end{array} \right) $ 
implies that $s_{11}t_{11}=1$ and
$s_{22}t_{22}=1$. Also, 
$TS=\left( \begin{array}{cc}
 1 & 0 \\
 0 & 1
\end{array} \right) $ 
implies that $t_{11}s_{11}=1$ and
$t_{22}s_{22}=1$. Thus, $t_{11}$ and $t_{22}$ have inverses $s_{11}$ and
$s_{22}$ in $B.$ By hypothesis, $s_{11},s_{22}\in A$. Also, \
$t_{11}s_{12}+t_{12}s_{22}=0$.  Therefore $s_{12}=-s_{11}t_{12}s_{22}\in A$.  \ This proves $S\in U_{2}(A)$.  

It is known that $T\in U_{n}(A)$ can have an inverse $S\in M_{n}(A)$
with $\ S\notin U_{n}(A)$; see for example the author's paper \cite{[B3]}.

These remarks lead to the question: \ If $A$ is inverse closed in $B$,  
then is $U_{n}(A)$ inverse closed in $M_{n}(B)$? Now we prove that this
question has an affirmative answer. In fact, we do the result for a
subalgebra larger than $U_{n}(A)$. We define this algebra now. 

Let $\mathit{I}(A)$ denote the largest inessential ideal of $A;$ see
[1, Fact~3.1]. To illustrate this concept in the case of
operators, let $B(X)$ and $K(X)$ denote the algebra of all bounded
linear operators on a Banach space $X$, and the space of all compact
operators on $X$, respectively. Then $\mathit{I}(B(X))$ is the closed
ideal of inessential operators on $X$ as originally defined by Kleinecke
\cite{[Kl]}. It is always true that $K(X)\subseteq \mathit{I}(B(X))$.  The key
fact here for our purposes is that the spectral theory of an element
$a\in \mathit{I}(A)$ is exactly like that of a compact operator. In
particular, for $a\in \mathit{I}(A)$,   $\sigma (a;A)$ is a finite set or
a sequence converging to 0; see [1, Remark~2.6]. This implies that
$a\in \overline{\hbox{Inv}(A)}$.  

\setcounter{definit}{4}
\begin{definit}$\left.\right.$\vspace{.3pc}

\noindent {\rm The matrix $T=(t_{jk})\in M_{n}(A)$ is in $U_{n}^{e}(A)$ if
$t_{jk}\in \mathit{I}(A)$ whenever $j>k$. The matrix $T=(t_{jk})\in
M_{n}(A)$ is in $L_{n}^{e}(A)$ if $t_{jk}\in \mathit{I}(A)$ whenever
$j<k$.}  
\end{definit}

\setcounter{theor}{5}
\begin{theor}[\!]
Assume that $A$ is inverse closed in $B$,  and $A$ is continuously
embedded in$\;B$. If $T\in $ $U_{n}^{e}(A)$ {\rm (}or $T\in L_{n}^{e}(A))$ and
$T^{-1}\in M_{n}(B)${\rm ,} then $T^{-1}\in M_{n}(A)$.  
\end{theor}

\begin{proof}
First assume that $n=2$.  Assume that $T\in U_{2}^{e}(A)$ has
$T^{-1}\in M_{2}(B)$.  By assumption, $t_{21}\in \mathit{I}(A)$,   so
$t_{21}$ is a limit of invertible elements in $A$. Thus Proposition~4
applies and so $T^{-1}\in M_{2}(A)$.  

The proof proceeds by induction. Assume that the result holds for $n$. 
Assume as a first case that: \ $T\in $ $M_{n+1}(A)$, $T^{-1}\in
M_{n+1}(B)$, $t_{jk}\in \mathit{I}(A)$ \ {\it provided that} $j>k$
{\it and} $ j\geq 3$, and $t_{11}$ $\in {\hbox{Inv}(A)}$.  \ We may assume
that $t_{11}=1$.  

Define $W$ and $V$ in $M_{n+1}(A)$ by:
\begin{alignat*}{3}
&w_{kk}&&=1,1\leq k\leq n+1;\quad w_{1k}=-t_{1k},\;2\leq k\leq
n+1;\\
&&&\qquad w_{jk}=0,\quad \hbox{otherwise;}\\
&v_{kk}&&=1,1\leq k\leq n+1; \quad v_{j1}=-t_{j1},\;2\leq k\leq
n+1;\\
&&&\qquad v_{jk}=0,\quad \hbox{otherwise.}
\end{alignat*}
It is easy to check that both $W$ and $V$ are invertible in
$M_{n+1}(A)$. Set $S=VTW$.  By direct computation
\begin{equation*}
s_{11}=1;\quad s_{1k}=0,\;2\leq k\leq n+1;\quad s_{j1}=0,\;2\leq j\leq n+1.
\end{equation*}
Define $\widetilde{S}\in M_{n}(A)$ by deleting the first row and the
first column of $S$.  A direct matrix computation verifies that
$\widetilde{S}\in U_{n}^{e}(A)$. Since $S^{-1}\in M_{n+1}(B)$,  
and noting the form of the first row and the first column of $S$
displayed above, we have\ $\widetilde{S}^{-1}\in M_{n}(B)$.  Then by
the induction hypothesis, $\widetilde{S}^{-1}\in $ $M_{n}(A)$.  It
follows that $S^{-1}\in M_{n+1}(A)$.  

In the remaining case, $t_{11}$ $\notin \hbox{Inv}(A)$.  Interchange rows 1 and
2 of $T$, so that $t_{21}$ is in the $(1,1)$ position. Call this new
matrix $R$, and note that $R\thickapprox T$. Since $r_{11}=t_{21}\in
\mathit{I}(A)$, we can choose $\lambda _{m}\rightarrow 0$ such that
$\lambda _{m}+r_{11}\in {\hbox{Inv}(A)}$ for all $m$. Let $R_{m}$ be the matrix
with $\lambda _{m}+r_{11}$ in the $(1,1)$-position, and with $r_{jk}$ in
the $(j,k)$-position otherwise. Since $R$ is invertible in $M_{n+1}(B)$
and $R_{m}\rightarrow R$ in norm, we may assume that $R_{m}^{-1}\in
M_{n+1}(B)$ for all $m$.  Now each $R_{m}$ satisfies the hypotheses of
the first case considered above, and so $R_{m}^{-1}\in $ $M_{n+1}(A)$
for all $m$.  Thus, $R^{-1}\in M_{n+1}(A)$ by applying Facts~1 and
2.\hfill \ba
\end{proof}

Now we look at a concrete situation where Theorem 6 applies. Let $(Y,\mu
)$ be a $\sigma $-finite measure space. Let $K(x,y)$ be a kernel with
the property that 
\begin{equation*}
|||\;K\;|||_{\infty }\equiv ess\sup_{x\in Y}\int_{Y}|K(x,y)|\;\hbox{d}\mu
(y)<\infty.
\end{equation*}
For such a kernel $K$, define the integral operator $T_{K}:L^{\infty
}\rightarrow L^{\infty }$ \ by
\begin{equation*}
T_{K}(f)(x)\equiv \int_{Y}K(x,y)f(y)\;\hbox{d}\mu (y),f \in L^{\infty }.
\end{equation*}
Then $T_{K}\in B(L^{\infty })$.  The set of all such integral operators
form a subalgebra of $B(L^{\infty })$ which is called the algebra of
Hille--Tamarkin operators on $L^{\infty };$ see [5, 11.5]. We use the
notation $H_{\infty }$ for this algebra. $H_{\infty }$ is a Banach
algebra with respect to the norm \ $|||K|||_{\infty }$.  This space
of integral operators is an important class of operators which contains
many interesting examples. 

We assume that $H_{\infty }$ does not contain the identity operator (the
usual situation). Let $H_{\infty }^{1}$ be the algebra $H_{\infty }$
with the identity operator on $L^{\infty }$ adjoined. It follows from
results in [5, 11.5] that $H_{\infty }^{1}$ is inverse closed in
$B(L^{\infty })$.  Let $ J=\{T\in H_{\infty }^{1}:T$ is a compact
operator on $L^{\infty }\}$.  Certainly, $J\subseteq
\mathit{I}(B(L^{\infty }))$.  Theorem 6 applies to any operator $\
T=(t_{jk})\in M_{n}(H_{\infty }^{1})$ with the property that $ t_{jk}\in
J$ \ whenever \ $j>k$.  

\section{Results when {\bi A} is a symmetric Banach *-algebra}

The following fact, which we use repeatedly, is proved in [10,
Theorem~9.8.4 and the preceding remarks, p. 1011]. 

\setcounter{fact}{6}
\begin{fact}
{\rm When $A$ is a symmetric Banach *-algebra, then $M_{n}(A)$ is
symmetric.} 
\end{fact}

\noindent The next proposition follows from a result of the author in \cite{[B1]} and 
Fact~7.

\setcounter{propo}{7}
\begin{propo}$\left.\right.$\vspace{.3pc}

\noindent When $A$ is a symmetric Banach *-algebra{\rm ,} $A$ is continuously embedded
in $B${\rm ,} and $A$ is closed in $B${\rm ,} then $M_{n}(A)$ is inverse closed in
$M_{n}(B)$ for all $n$.  
\end{propo}

\begin{proof}
By Fact 7, $M_{n}(A)$ is symmetric. Also, $M_{n}(A)$ is continuously
embedded in $M_{n}(B)$, and $M_{n}(A)$ is closed in $M_{n}(B)$. Then
the proposition follows from [2, Theorem]. \hfill \ba
\end{proof}
\setcounter{theor}{8}
\begin{theor}[\!]
Assume that $A$ is a symmetric Banach *-algebra{\rm ,} that $A$ is *-inverse
closed in $B${\rm ,} and that $A$ is continuously embedded in $\;B$. Then
$M_{n}(A)$ is symmetric and *-inverse closed in $M_{n}(B)$.  
\end{theor}

\begin{proof}
$M_{n}(A)$ is symmetric by Fact 7.

Assume that $n=2$. Suppose that $T=T^{\ast }\in M_{2}(A)$ and that $
T^{-1}\in M_{2}(B)$.  Write $T=(t_{jk})$,   and note that
$t_{11}=t_{11}^{\ast }$.  Since $A$ is symmetric, $t_{11\text{ \ }}$is
the limit of the sequence of invertible elements
$\{\frac{i}{n}+t_{11}\}$.  By Proposition 4, it follows that $T^{-1}\in
M_{2}(A)$.  

Now note that $M_{2}(M_{2^{n}}(A))$ can be naturally identified as a
Banach *-algebra with $M_{2^{n+1}}(A)$.  Therefore by the case where
$n=2$ proved above, and induction we have: 
\begin{equation*}
M_{2^{n}}(A)\ \hbox{is *-inverse closed in}\ M_{2^{n}}(B)\ \hbox{for}\ n\geq 1.
\end{equation*}

Assume that $2^{n-1}<m<2^{n}$ for some $n\geq 2$. \, For $T\!=\!(t_{jk})\in
M_{m}(A)$,   define $\widetilde{T}\!=\!(\widetilde{t}_{jk})$ $\in
M_{2^{n}}(A)$ by\vspace{.6pc} 
\begin{equation*}
\widetilde{t}_{jk}=t_{jk},1\leq j,k\leq m;\quad 
\widetilde{t}_{kk}=1,\;m<k\leq 2^{n};\quad \widetilde{t}_{jk}=0,\quad
\hbox{otherwise.}
\end{equation*}
It is straightforward to check that $T\rightarrow \widetilde{T}$ is a
unital *-algebra monomorphism of $M_{m}(A)$ into \ $M_{2^{n}}(A)$, and
that $T$ is invertible in $M_{m}(A)$ if and only if $\widetilde{T}$ is
invertible in $ M_{2^{n}}(A)$. The same definition as given above defines
a unital algebra monomorphism of $M_{m}(B)$ into \ $M_{2^{n}}(B)$, and
again, $T$ is invertible in $M_{m}(B)$ if and only if $\widetilde{T}$ is
invertible in $ M_{2^{n}}(B)$.  

Finally, if $T=T^{\ast }\in M_{m}(A)$ and $T^{-1}\in M_{m}(B)$, then $
\widetilde{T}=\widetilde{T}^{\ast }\in M_{2^{n}}(A)$ and that
$\widetilde{T}^{-1}\in $ $M_{2^{n}}(B)$.  As proved previously, this
implies $\widetilde{T}^{-1}\in M_{2^{n}}(A)$,   and so $T^{-1}\in
\break M_{m}(A)$.
\hfill \ba
\end{proof}\vspace{-.2pc}

The following corollary extends [3, Theorem~11] which is the case $n=1$. 

\pagebreak

\setcounter{coro}{9}
\begin{coro}$\left.\right.$\vspace{.3pc}

\noindent Assume that $A$ is a symmetric Banach *-algebra{\rm ,} and that $A$
is either *-inverse closed in $B$ and continuously embedded in$\;B${\rm ,} or
$A$ is \hbox{\rm SRP} in $B$. Also{\rm ,} assume that for some $M>0${\rm ,}   
\end{coro}
\begin{equation*}
\Vert a^{\ast }\Vert \leq M\Vert a\Vert  \quad \hbox{for all}\ a\in A.
\end{equation*}
{\it Then} $M_{n}(A)$ {\it is inverse closed in} $M_{n}(B)$ {\it for} $n\geq 1$.

\begin{proof}
Since $A$ is symmetric, $M_{n}(A)$ is symmetric by Fact 7. If $A$ is
*-inverse closed and continuously embedded in $B$,   then by Theorem 9,
$M_{n}(A)$ is *-inverse closed in $M_{n}(B)$.  If $A$ is SRP in $B$,
then by Fact 1, $M_{n}(A)$ is SRP in $M_{n}(B)$.  Also, it is easy to
check that for all $T\in M_{n}(A)$,   $\Vert T^{\ast }\Vert \leq M\Vert
T\Vert $. In both cases it follows from [3, Theorem~11] that $M_{n}(A)$ is
inverse closed in $M_{n}(B)$. \hfill \ba
\end{proof}



\begin{thebibliography}{9999}
\bibitem{[B0]} Barnes B, Murphy G, Smyth R and West T, Riesz and
Fredholm theory in Banach Algebras, {\it Pitman Research Notes in Math}, {\bf 67},
Boston, 1982

\bibitem{[B1]} Barnes B, A note on the invariance of spectrum for
symmetric Banach *-algebras, {\it Proc. Am. Math. Soc.} {\bf 126} (1998)
3545--3547

\bibitem{[B2]} Barnes B, Symmetric Banach *-algebras: invariance of
spectrum, {\it Studia Math.} {\bf 141} (2000) 251--261

\bibitem{[B3]} Barnes B, The spectral theory of upper triangular matrices
with entries in a Banach algebra, {\it Math.\ Nachrichten} {\bf 241}
(2002) 5--20

\bibitem{J} Jorgens K, Linear integral operators (Boston: Pitman)
(1982)

\bibitem{[K]} Krupnik N, Banach algebras with symbol and singular integral
operators (Basel: Birkhauser-Verlag) (1987)

\bibitem{[KM]} Krupnik N and  Markus M, On the inverse closedness of certain
Banach subalgebras, {\it Studies on Diff. Euat. and Math. Anal.},
'Stiintsa', Kishinev (1988) 93--99 (Russian)

\bibitem{[Kl]} Kleinecke D, Almost-finite, compact, and inessential operators,
{\it Proc. Am. Math. Soc.} {\bf 14} (1963) 863--868

\bibitem{[P1]} Palmer T, Banach algebras and the general theory of *-algebras
(Cambridge: Cambridge Univ. Press) (1994) vol. 1

\bibitem{[P2]} Palmer T, Banach algebras and the general theory of
*-algebras (Cambridge: Cambridge Univ. Press) (2001) vol. 2
\end{thebibliography}
\end{document}